
\baselineskip=14pt
\parskip=10pt

\magnification=\magstephalf

\def\P{{\cal P}}

\def\1{{\overline{1}}}
\def\2{{\overline{2}}}
\parindent=0pt
\overfullrule=0in

\def\frac#1#2{{#1 \over #2}}
\centerline
{\bf 
Enumerating Seating Arrangements that Obey Social Distancing
}
\bigskip
\centerline
{\it George SPAHN and Doron ZEILBERGER}
\centerline
\qquad \qquad \qquad  {\it In fond memory of Marko Petkovsek (1955-2023)}
\bigskip
{\bf Abstract}:  
We illustrate the power of symbolic computation and experimental mathematics by  investigating
maximal seating arrangements, either on a line, or in a rectangular auditorium with a fixed number of columns but an arbitrary number of rows, 
that obey any prescribed set of `social distancing' restrictions.
In addition to enumeration, we study the statistical distribution of the density, and give simulation algorithms for generating them.



{\bf Preface: How it all started} 

It all started when we came across the delightful article [PSZ] whose motivation had nothing to do with the now fashionable {\it social distancing}. Their starting point
was an $r \times s$  housing development with $rs$ building lots, each with room for a house. They wanted to enumerate the
number of ways (out of the total of $2^{rs}$ possibilities, including building nothing, and building on all the lots) for which

$\bullet$ No house is `blocked from the sun' (see below).

$\bullet$ You can't build any house on a currently empty lot without violating this condition.

They were not only interested in the number of such configurations, let's call it $T(r,s)$, but among those, the statistical
distribution of the number of houses, or equivalently, the density (the number of houses divided by $rs$), in other
words, the {\it generating function}, or {\it weight-enumerator}, of such maximal configurations according to the
weight $z^{NnmberOfHouses}$.

It turns out that from an abstract point of view, this is nothing but enumerating {\bf maximal} (w.r.t. the number of ones) $r \times s$ $0-1$ matrices that avoid the pattern
$$
\matrix{1 &  1 & 1 \cr
          &  1 & } \quad.
$$

In other words, $0-1$ matrices $(a_{ij})$ ($1 \leq i \leq r$, $1 \leq j \leq s$) where it is forbidden to have a location $(i_0,j_0)$, where
$$
a_{i_0,j_0-1}=1 \quad, \quad  a_{i_0,j_0}=1 \quad, \quad a_{i_0,j_0+1}=1 \quad, \quad a_{i_0+1,j_0}=1 \quad ,
$$
and changing any of the entries that are currently $0$ to $1$, would create this undesirable pattern.

But why just this particular pattern? The question makes sense for enumerating such $0-1$ matrices avoiding {\it any} pattern, and in fact, any {\it set of patterns}.

In particular this question has immediate relevance to {\it social distancing}. We have a rectangular classroom with $r$ rows, each with $s$ seats, 
and you can't have two students sitting next to each other in the same row, or having anyone immediately in front of you, or behind you (i.e. you can't have any two students sitting next to each other
in the same column). In addition
none of the currently empty seats can accommodate a newcomer without breaking this restriction. This  is equivalent to the problem of weight-enumerating $r \times s$
{\bf maximal} (with respect to the number of ones) $0-1$ matrices avoiding the patterns
$$
\matrix{ 1 & 1} \quad \quad \quad \quad , \quad \quad \quad  \matrix{1 \cr 1} \quad .
$$

The problem is already interesting for one row (equivalently one column), and one can impose  lots of possible restrictions.
Going back to two dimensions, one can think of prohibiting the pattern
$$
\matrix{ 1 &1 \cr 1 &1 } \quad .
$$
and the possibilities are endless.

The famous problem of {\it maximal non-attacking kings} can be formulated as counting (and weight-counting according to the number of $1$s)
maximal $r \times s$ $0-1$ matrices avoiding the patterns
$$
\matrix{1 & 1} \quad \quad , \quad \quad
\matrix{1 \cr 1} \quad \quad , \quad
\matrix{1 & \cr  & 1} \quad \quad , \quad \quad
\matrix{ & 1\cr 1 & } \quad.
$$

We will probably never know the exact number of doing it for a $100 \times 100$ board, and in general, for any fixed pattern, or set of patterns, the problem of
enumerating such maximal $0-1$ {\bf square} matrices, seems intractable. But thanks to the {\bf transfer matrix} method one can efficiently find
explicit bi-variate generating functions for such enumeration problems with a {\bf fixed} number of columns, $s$ (not too big, of course), but arbitrary number of rows $r$ (or vice versa).
So the problem of enumerating (and weight-enumerating) maximal configurations of non-attacking kings, say on a $1000 \times 5$ board can be found exactly.

To see their total number, and the number of these maximal placements with $668$ kings (the smallest possible number),
see the following output file:

{\tt https://sites.math.rutgers.edu/\~{}zeilberg/tokhniot/oSDk5Story.txt} \quad .

We will now describe how to do it.

{\bf A Finite State Machine for Maximal Seating}

Suppose that we have $r$ rows of seats. We wish to compute the sequence giving the number of maximal seating assignments as a function of the number of columns. As input we receive $r$, the number of rows, and $\P$, the set of 
patterns describing the shapes to be avoided.

There are 2 ways a seating arrangement can fail.

$\bullet$  It may fail to meet the social distancing requirement. That is there may be some people sitting in a place prohibited by one of the patterns. Call this a failure of type 1.

$\bullet$ It may fail to be maximal. That is there may be an empty seat which could be filled without violating any of the patterns. Call this a failure of type 2.

We note that both of these problems are {\it local}. We don't need to know the entire seating arrangement to determine with certainty that a possible assignment is not valid. This allows us to create a finite state machine with relative ease. We treat a 
seating arrangement as a sequence of columns. The possible columns will serve as inputs to our state machine. Here, we represent a column as a binary string of 0s and 1s of length $r$. 
The 1s correspond to occupied seats, and the 0s correspond to empty seats. There are thus $2^r$ possible symbols that our state machine must be able to process.

As we read in columns, we must determine whether the constraints have been satisfied {\it so far}. In order to detect whether the proposed next column causes a violation, 
we need to store some data about the preceding columns. In this case, we just need to store the contents  of some number of previous columns.

{\bf Definition:} For a given pattern, $p$, define it's width, $w(p)$, to be the largest $x$-value of any of its lattice coordinates minus the smallest $x$ value. Define 
the width of our state machine, $W(r,\P)$ or just $W$, to be the greatest width of all the pattern in $\P$.

To check whether the proposed next column causes a violation of type 1, we need to have stored the previous $W$ columns. 
If we have that information stored, we can just loop through all the patterns and check that each new 1 does not add a violation of that pattern shape. 

To check whether the proposed next column causes a violation of type 2, we must be a bit more clever. For each 0 that is added, there must be a group of nearby 1's that ensures it cannot be changed to a 1 without causing a violation. 
A 0 that has an appropriately shaped group of nearby 1s is called a satisfied 0. Since we do not know the contents of columns we might read in the future, we are unable to determine whether the 0s 
in the current column are satisfied. Sure, they may not be satisfied yet, but we should not consider them to be 
violations yet. If the current column is column $j$, then we check that the 0s in column $j-W$ are satisfied. 

In total, we store the contents of the previous $2W$ columns in the state machine. Thus for column $j-W$ we know exactly the contents of the previous $W$ columns and the next $W$ columns so we can determine with certainty whether the 0s are satisfied.

Since we require a state for each of the possible $2W$ previous columns, this gives a total of a $2^{2Wr}$ possible states.

Now that we have constructed our set of states, it is simple to determine the transitions. For each state, we have a transition for each possible next column. If the transition causes no violations, it is valid, otherwise the state machine should immediately REJECT. Once a violation has been detected we do not need to look at any more columns.

When we are done reading our input, we must do a little more work to determine whether to ACCEPT or REJECT. Specifically, the 0s in the most recent $W-1$ columns have not been checked yet, so we must check that each is satisfied before accepting. 

We also add some special initialization states for when the previous $2W$ columns do not have values yet. The details of this are implementation specific and omitted.

{\bf Constructing the Transfer Matrix}

Once the finite state machine has been constructed, we seek to count the number of input strings of length $n$ that cause the machine to ACCEPT. To do this we use the transfer matrix method. Each state corresponds to a row and column of the matrix, $M$. 
Then the $i,j$ entry $M[i][j]$ is 1 if there is a valid transition from state $i$ to state $j$, and 0 otherwise. 

We can count the number of paths from state $i$ to state $j$ of length $s$ by computing $M^s [i][j]$. Thus the number of maximal seating arrangements with $s$ columns is given by an entry of the matrix $M^s$. We just set $i$ to 
be the initial state and $j$ to be the final state of our state machine.

{\bf Density of Seating arrangements}

We can modify the transition matrix by replacing each $1$ with a power of $z$. 
For a transition corresponding to adding a column with $t$ 1s, the corresponding matrix entry will be $z^t$. Then $M^s[i][j]$ is a polynomial in $z$. The coefficient of $z^k$ is the number of seating arrangements with a total of $k$ seats occupied.

{\bf Getting the bi-variate generating function}

The previous section explained how to compute the transfer matrix, let's call it $M_{\P,r}(z)$, for {\it any} set of patterns $\P$, and any {\bf fixed} number of rows $r$. 
Altogether it had, say, $d$, states, including the two special states, the initial, and the final, that  our program labeled $1$ and $d$, respectively.

We are interested, for an {\it arbitrary} number of columns, $s$,
in the weight-enumerator, let's call it $W_{\P,r,s}(z)$, according to the weight $z^{NumberOfOnes}$ of all maximal $r \times s$ $0-1$ matrices avoiding $\P$. Define the bi-variate generating function
$$
f_{\P,r}(z,x) := \sum_{s=0}^{\infty} \, W_{\P,r,s}(z) x^s \quad.
$$
Then we have
$$
f_{\P,r}(z,x) = (I\,-\, x\,M_{\P,r}(z))^{-1}[1,d] \quad.
$$
What's nice about computer algebra is that it can all be streamlined.

Using this generating function, we can automatically compute the {\it limiting average density} (among maximal seatings) for a fixed $r$ and arbitrary $\P$, as $s$ goes to infinity, namely
$$
 \frac{1}{r}  \lim_{s \rightarrow \infty} \,\frac{ W_{\P,r,s}'(1)}{s\, W_{\P,r,s}(1)} \quad.
$$

We omit the details, since this is standard residue calculus, and the interested reader can look at the Maple source code of procedure {\tt AsyAv(f,x,z)} in the Maple package {\tt SD.txt}, accompanying this paper.

{\bf Random Sequential Adsorption}

The polynomial $W_{\P,r,s}(z)$ is the weight enumerator of maximal $0-1$ $r \times s$ matrices that avoid the patterns in $\P$. Hence the coefficient of $z^m$ in that polynomial is the
exact number of such matrices with exactly $m$ ones (or equivalently the number of maximal seatings in an $r$ by $s$ classroom obeying the restrictions in $\P$ and having exactly $m$
occupied chairs). It follows that the {\it probability generating function} (for the uniform distribution), let's call it $V_{\P,r,s}(z)$, is:
$$
V_{\P,r,s}(z)=\frac{W_{\P,r,s}(z)}{W_{\P,r,s}(1)} \quad,
$$
and the expected number of occupied seats is  $\frac{d}{dz} V_{\P,r,s}(z)|_{z=1}= V_{\P,r,s}'(1)$.

There is a quick way to generate random maximal seatings, described in the special case of $T$-avoiding seatings in [PSZ], but that makes sense in general, called
{\it Random Sequential Adsorption}, that they abbreviated to RSA, but since for us RSA means `Rivest-Shamir-Adleman', we will refrain from using this abbreviation.

In that random-generation process, one picks, {\it uniformly at random}, a permutation of the $rs$ initial empty seats, representing $rs$ people each having their
favorite seat, all distinct from each other.
When someone enters the classroom they attempt to occupy their favorite seat, and do so if they do not cause a violation.
But if sitting there does create a violation, they are not allowed to sit elsewhere and must leave the room (and cut the class).

This process is very easy  to simulate, and it generates a random maximal seating. Alas, it is not the uniform one. Thanks to our closed-form expressions for
the probability distribution for $r \times s$ maximal seatings with a fixed (not too big) $r$ but arbitrarily large $s$, we were able to compare
notes, and estimate (via simulations) how close is random sequential adsorption to the uniform distribution.

{\bf The Maple package SD.txt}

This article is accompanied by a Maple package {\tt SD.txt}, available from 

{\tt https://sites.math.rutgers.edu/\~{}zeilberg/tokhniot/SD.txt} \quad .

The front of this article 

{\tt https://sites.math.rutgers.edu/\~{}zeilberg/mamarim/mamarimhtml/social.html} \quad,

contains numerous sample output files, a  selection of which we will describe next.

{\bf Output}

{\bf Explicit Expressions for generating functions in One-Dimensional Seatings: Avoiding $b$ consecutive occupied seats}

Let $A_b(m,n)$ be the number of maximal $0-1$ vectors of length $n$ with $m$ ones (equivalently: sum $m$), avoiding $b$ consecutive $1$'s, and let
$$
f_b(z,x):=\sum_{n=0}^{\infty} \, \left ( \sum_{m=0}^{n} A_b(m,n)\, z^m \right ) \,x^n \quad,
$$
then we have
$$
f_2(z,x)=-\frac{x^{2} z +x z +1}{x^{3} z +x^{2} z -1} \quad,
$$
and the limiting average density is $.4114955887\dots$ \quad .

$$
f_3(z,x)=-\frac{x^{5} z^{3}-x^{3} z^{2}-x^{2} z^{2}+x^{2} z -x z -1}{x^{6} z^{3}-x^{4} z^{2}-x^{3} z^{2}-x^{2} z +1} \quad,
$$
and the limiting average density is $0.5772029462\dots$ \quad.

$$
f_4(z,x)=
-\frac{x^{9} z^{6}+x^{7} z^{5}-x^{5} z^{4}+x^{5} z^{3}-2 x^{4} z^{3}-x^{3} z^{3}+x^{3} z^{2}-x^{2} z^{2}-x z -1}{x^{10} z^{6}+x^{8} z^{5}-x^{6} z^{4}-2 x^{5} z^{3}-x^{4} z^{3}-x^{3} z^{2}+1}
$$
and the limiting average density is $0.6686427921\dots$ \quad.

$$
f_5(z,x)=
$$
$$
-\frac{x^{14} z^{10}-x^{11} z^{8}-2 x^{9} z^{7}+x^{9} z^{6}-2 x^{8} z^{6}+x^{6} z^{5}-x^{6} z^{4}+2 x^{5} z^{4}+x^{4} z^{4}-2 x^{4} z^{3}+x^{3} z^{3}-x^{3} z^{2}+x^{2} z^{2}+x z +1}{x^{15} z^{10}-x^{12} z^{8}-2 x^{10} z^{7}-2 x^{9} z^{6}+x^{7} z^{5}+2 x^{6} z^{4}+x^{5} z^{4}+x^{4} z^{3}+x^{3} z^{2}-1}
$$
and the limiting average density is $0.7269949175\dots$ \quad.

For $f_6(z,x),f_7(z,x),f_8(z,x)$, see the output file

{\tt https://sites.math.rutgers.edu/\~{}zeilberg/tokhniot/oSD1.txt} \quad .

The limiting average densities are $0.7675902978$, $0.7975140257$, and $0.8205096203$, respectively.

\vfill\eject

{\bf Explicit Expressions for generating functions in One-Dimensional Seatings: Stricter Social Distancing Restrictions}

Suppose that, in a row of seats of length $n$,  each person must have at least $b$ empty seats on either side, and
let $C_b(m,n)$ be the number of such maximal $0-1$ vectors of length $n$ and $m$ $1$s, and define the bi-variate generating function:
$$
g_b(z,x):=\sum_{n=0}^{\infty} \, \left ( \sum_{m=0}^{n} C_b(m,n)\, z^m \right ) \,x^n \quad .
$$
To our pleasant surprise, there is an {\bf explicit} expression not only in terms of $z$ and $x$ but also for symbolic $b$. We leave it as challenge to the reader
to prove it for arbitrary $b$.

{\bf Conjecture}: For $b \geq 1$, 
$$
g_b(z,x) \, = \,
{\frac {{x}^{2\,b+2}z-{x}^{b+1}z-{x}^{b+2}z+xz+ \left( -1+x \right) ^{2}}{ \left( -1+x \right)  \left( -{x}^{2\,b+2}z+{x}^{b+1}z+x-1 \right) }} \quad .
$$

We confirmed it for $b \leq 8$. See the output file

{\tt https://sites.math.rutgers.edu/\~{}zeilberg/tokhniot/oSD2.txt} \quad .

{\bf Avoiding Dimers in $3 \times s$,  $4 \times s$, and  $5 \times s$, $0-1$ matrices}

To see explicit expressions, limiting average densities, and comparison with Random Sequential Adsorption simulation, for the $3 \times s$ case, see the output file

{\tt https://sites.math.rutgers.edu/\~{}zeilberg/tokhniot/oSDd3.txt} .

In particular, the bi-variate generating function for  $3 \times s$ maximal dimer-avoiding $0-1$ matrices is:
$$
-\frac{\left(2 x^{5} z^{4}+2 x^{3} z^{3}+2 x^{2} z^{3}-6 x^{2} z^{2}-x \,z^{2}-x z -z -1\right) x z}{x^{5} z^{4}+2 x^{4} z^{4}-x^{4} z^{3}+x^{3} z^{4}-4 x^{3} z^{3}-x^{2} z^{3}-x z +1} \quad .
$$

For similar information for the $4 \times s$ case, see the output file

{\tt https://sites.math.rutgers.edu/\~{}zeilberg/tokhniot/oSDd4.txt} ,

in particular, the bi-variate generating function for  $4 \times s$ maximal dimer-avoiding $0-1$ matrices: is,
$$
\frac{\left(x^{6} z^{6}-x^{5} z^{6}+x^{5} z^{5}-2 x^{5} z^{4}-3 x^{4} z^{4}+2 x^{3} z^{4}-7 x^{3} z^{3}+2 x^{3} z^{2}-4 x^{2} z^{3}+7 x^{2} z^{2}-x \,z^{2}+4 x z +3\right) x \,z^{2}}{x^{6} z^{6}+x^{5} z^{6}+x^{5} z^{5}+2 x^{4} z^{5}-x^{4} z^{4}+2 x^{3} z^{5}-4 x^{3} z^{4}-x^{3} z^{3}-2 x^{2} z^{3}-x^{2} z^{2}-x \,z^{2}+1}
$$

For similar information for the $5 \times s$ case, see the output file:

{\tt https://sites.math.rutgers.edu/\~{}zeilberg/tokhniot/oSDd5.txt} .

{\bf Maximal Non-Attacking Kings}

The generating function for maximal configurations of non-attacking kings on a $3 \times s$ chessboard is
$$
-\frac{\left(x^{5} z^{3}+x^{5} z^{2}-x^{3} z^{3}+x^{3} z +2 x^{2} z^{2}+x^{2} z +x \,z^{2}-x^{2}-3 x z -2 x -z -1\right) x z}{x^{6} z^{4}+x^{6} z^{3}-x^{5} z^{4}-x^{5} z^{3}+x^{4} z^{3}+x^{4} z^{2}+x^{3} z^{3}-x^{3} z^{2}-x^{3} z -x^{2} z^{2}-x^{2} z -x z +1} \quad .
$$
  
For more details see:
{\tt https://sites.math.rutgers.edu/\~{}zeilberg/tokhniot/oSDk3.txt} \quad .

The generating function for maximal configurations of non-attacking kings on a $4 \times s$ chessboard is
$$
-\frac{\left(6 x^{6} z^{3}+9 x^{5} z^{2}-6 x^{4} z^{3}+3 x^{4} z^{2}-3 x^{3} z^{2}+3 x^{3} z +3 x^{2} z^{2}+2 x^{2} z -3 x^{2}+3 x z -12 x -3\right) x \,z^{2}}{6 x^{7} z^{5}-6 x^{6} z^{5}+9 x^{6} z^{4}-6 x^{5} z^{4}+3 x^{4} z^{4}+x^{4} z^{3}+3 x^{3} z^{3}-6 x^{3} z^{2}-4 x^{2} z^{2}-x z +1} \quad .
$$

For more details see:
{\tt https://sites.math.rutgers.edu/\~{}zeilberg/tokhniot/oSDk4.txt} \quad .

For the bi-variate generating function of maximal configurations of non-attacking kings on a $5 \times s$ chessboard  and more details, and results of simulations, see:

{\tt https://sites.math.rutgers.edu/\~{}zeilberg/tokhniot/oSDk5.txt} \quad .

See also the output files

{\tt https://sites.math.rutgers.edu/\~{}zeilberg/tokhniot/oSD22C3.txt} \quad,

{\tt https://sites.math.rutgers.edu/\~{}zeilberg/tokhniot/oSD22C4.txt} \quad,

for avoiding a $2 \times 2$ block of $ones$.

Readers are welcome to generate more data, and experiment with other patterns (and set of patterns) to their heart's content, using the Maple package {\tt SD.txt}. Enjoy!

{\bf The Uniform   vs. the   Random Sequential Adsorption Distributions}

It seems much harder to get the {\it exact} probability generating functions for the density (equivalently, the number of $1$s)  by using  Random Sequential Adsorption
in general, and we have to resort to (very reliable!) simulations.

But for the very special case of seatings on a line with $n$ seats, and avoiding the single pattern $11$, i.e. non adjacent seatings, we can get exact, closed-form expression
for {\bf both} the uniform and the random sequential adsorption process.

Let $g_n(z)$ be the enumeration generating function, according to `number of occupied seats' under the uniform distribution, then as we saw above
$$
\sum_{n=0}^{\infty} \, g_n(z)\, x^n =  -\frac{x^{2} z +x z +1}{x^{3} z +x^{2} z -1} \quad .
$$
Equivalently the sequence $\{g_n(z)\}$  satisfies the {\bf linear} recurrence
$$
g_n(z)= z (g_{n-2}(z)+ g_{n-3}(z)) \quad,
$$
subject to the initial conditions
$$
g_0(z)=1 \quad, \quad
g_1(z)=z \quad, \quad
g_2(z)=2z \quad.
$$

To turn it into a {\it probability generating function}, we have to divide by $g_n(1)$ and define
$$
G_n(z):=\frac{g_n(z)}{g_n(1)} \quad .
$$

Recall that in random sequential adsorption the input is a permutation $\pi$. If $\pi(1)=i$ then the persons that like locations $i-1$ and $i+1$ would not be able to sit, since $i$ is already occupied. There
are $(n-1)(n-2)$ ways where they can be placed. The remaining students independently seat in locations $1\leq x\leq i-2$ and  $i+2\leq x\leq n$ and it follows that
the probability generating function, let's call it $F_n(z)$ (after dividing by $n!$),  satisfies the {\bf non-linear} recurrence
$$
F_n(z) = \frac{z}{n} \left ( 2\, F_{n-2}(z) + \sum_{i=2}^{n-1} F_{i-2}(z) F_{n-i-1}(z) \right ) \quad , \quad n\geq 2 \quad ; \quad F_0(z)=1 \quad , \quad F_1(z)=z \quad .
$$

To get the expected number of occupied chairs under the uniform and the sequential adsorption distributions, we compute $G'_n(1)$ and $F'_n(1)$ respectively, or more informatively, the {\it densities} 
$G'_n(1)/n$ and $F'_n(1)/n$. They both converge to limits. The former can be found exactly. We have (as already stated above, with less precision, the number is a certain algebraic number)
$$
\lim_{n \rightarrow \infty} \frac{G'_n(1)}{n} \, =\,
0.411495588662645763381900381335531940800608649354765817635803356939840\dots
$$

The other limit we estimated numerically. we have (by taking large $n$)
$$
\lim_{n \rightarrow \infty} \frac{F'_n(1)}{n} \, \approx \, 0.4323\dots  \quad,
$$
agreeing with the exact value $(1-e^{-2})/2$ stated in Steven Finch's wonderful book [F] (section 5.3.1) and due to several people (see [F] for references). 

Hence the average occupation densities, for large $n$,  of the random sequential adsorption process is approximately $1.052$ times those of the 
uniform distribution. Of course, the maximal density, obtained by the
regular spacing $101010\dots$ is $0.5$.

\vfill\eject

{\bf Conclusion}

In this case study, we illustrated the power of {\it symbolic computation} and {\it experimental mathematics}, both so dear to Marko Petkovsek,
to get insight about the statistics of maximal seating arrangements avoiding any set of social distancing restrictions.

{\bf References}

[F] Steven R. Finch, {\it Mathematical Constants} (Encyclopedia of mathematics and its applications, v. 94) , Cambridge University Press, 2003.

[PSZ] Mate Puljiz, Stjepan Sebek, and Josip Zubrinic, {\it Packing density of combinatorial settlement planning models},
Amer. Math. Monthly {\bf 130}(10) (Dec. 2023).

[V] Vince Vatter, {\it Social Distancing, Primes, and Perrin Numbers}, , Math Horizons, {\bf 29}(1) (2022). \hfill\break
{\tt https://sites.math.rutgers.edu/\~{}zeilberg/akherim/vatter23.pdf} \quad .

\bigskip
\hrule
\bigskip
George Spahn and Doron Zeilberger, Department of Mathematics, Rutgers University (New Brunswick), Hill Center-Busch Campus, 110 Frelinghuysen
Rd., Piscataway, NJ 08854-8019, USA. \hfill\break
Email: {\tt  gs828 at math dot rutgers dot edu} \quad, \quad {\tt DoronZeil] at gmail dot com}   \quad .

Written: {\bf Jan. 22, 2024}. 

\end